\documentclass[11pt, eqno, twoside]{article}
\usepackage{amsfonts, amsmath, amssymb, amscd}
\usepackage{latexsym}
\usepackage{graphicx}
\usepackage{amsthm}

\allowdisplaybreaks[4]

\theoremstyle{plain}
\newtheorem{theorem}{Theorem}[section]
\newtheorem{lemma}[theorem]{Lemma}

\theoremstyle{definition}

\newtheorem{remark}[theorem]{Remark}

\numberwithin{equation}{section}

\title{
$L^2$-harmonic $p$-forms on submanifolds with finite total curvature
\thanks{This work is supported by NSF of Anhui Provincia Education
Department (No. KJ2017A341) and Talent Project of Fuyang Normal
University (No. RCXM201714)}}
\author {Jundong Zhou$^{a,\,b}$\\ \thanks{Email address:  zhoujundong109@163.com }
 $^a${\small School of Mathematics and Statistics, Fuyang Normal University, Fuyang, }\\
 {\small Anhui 236041, People's Republic of China}\\
 $^b${\small School of Mathematical Sciences,  University of Science and Technology of China, Hefei }\\
 {\small Anhui 230026, People's Republic of China}}

\date{}
\begin{document}
\maketitle

\begin{abstract}
\noindent Let $H^p(L^2(M))$ be the space of all $L^2$-harmonic
$p$-forms $(2\leq p\leq n-2)$ on complete submanifolds $M$ with flat
normal bundle in spheres. In this paper,  we first show that
$H^p(L^2(M))$ is trivial if the total curvature of $M$ is less than
a positive constant depending only on $n$. Second, we show that the
dimension of $H^p(L^2(M))$ is finite if the total curvature of $M$
is finite. The vanishing theorem is a generalized version of
Gan-Zhu-Fang theorem and the finiteness theorem is an extension of
Zhu-Fang theorem.
\end{abstract}



\section{Introduction}
$L^2$-harmonic forms on submanifolds  in various ambient spaces had
been studied extensively during past few years. Many results
demonstrated the fact that there is a close relation between the
topology of the submanifold and the total curvature by using theory
of $L^2$-harmonic forms. In \cite{5,6}, It was showed that a
complete minimal hypersurface $M^m$  $(m\geq3)$ with sufficient
small the total scalar curvature in $R^{n+1}$ has only one end. In
2008, Seo \cite{7} improved the upper bound of the total scalar
curvature which was given by Ni \cite{6}. Later Seo \cite{8} proved
that if an $n$-dimensional complete minimal submanifold $M$ in
hyperbolic space has sufficiently small total scalar curvature, then
$M$ has only one end. It is well-known that Euclidean space and
hyperbolic space are space forms all. Fu and Xu \cite{9} studied
$L^2$-harmonic 1-forms on complete submanifolds in space forms and
proved that a complete submanifold $M^n$ $(n\geq3)$  with finite
total curvature and some conditions on mean curvature must have
finitely many ends. Furthermore, Cavalcante, Mirandola and
Vit\'{o}rio \cite{10} obtained that if $M^n$ $(n\geq3)$ is a
complete noncompact submanifold in Cartan-Hadamard manifold with
finite total curvature and the first eigenvalue of the Laplacian of
$M^n$ is bounded from below by a suitable constant, then the space
of the $L^2$-harmonic $1$-forms on $M^n$ has finite dimension. Zhu
and Fang \cite{11} investigated complete noncompact submanifolds in
a sphere and obtained a result which was an improvement of Fu-Xu
theorem on submanifolds in spheres. Meanwhile, Zhu-Fang result was a
 generalized version of Cavalcante, Mirandola and Vit\'{o}rio's result
on submanifolds in Hadamard manifolds. The following theorem $A$ is
Zhu-Fang result.
 \vskip
5pt\noindent{\bf {Theorem A }}( [11] ) {\it Let $M^n (n \geq 3)$ be
an $n$-dimensional complete noncompact oriented manifold
isometrically immersed in an $(n+ p)-$dimensional sphere $S^{n+p}$.
If the total curvature is finite, then the dimension of
$H^1(L^2(M))$ is finite and there are finitely many non-parabolic
ends on M}.

In 2015, Lin \cite{14} studied $L^2$-harmonic $p$-forms on complete
submanifolds in Euclidean space and proved that if a complete
submanifold $M^n (n\geq3)$ with flat normal bundle in $R^{n+p}$ has
sufficient small the total curvature, then the space of the
 $L^2$-harmonic $p$-forms on $M^n$ is  trivial.
Recently, Gan, Zhu and Fang [20] studied $L^2$-harmonic $2$-forms on
complete noncompact minimal hypersurface in spheres and proved the
following result.

 \vskip 5pt\noindent{\bf {Theorem B}}( [20] )
{\it  Let $M^n (n \geq 3) $ be an $n$-dimensional complete
noncompact minimal hypersurface isometrically immersed in an $(n +
1)$-dimensional sphere $S^{n+1}$. There exists a positive constant
$\delta(n)$ depending only on $n$ such that if the total curvature
is less than $\delta(n)$, then the second space of reduced $L^2$
cohomology of M is trivial }.

Inspired  by Li-Wang work \cite{4} and the above results, in this
paper, we study the space of $L^2$-harmonic $p$-forms on submanifold
in spheres and prove the following vanishing and finiteness
theorems.

\vskip 5pt\noindent{\bf {Theorem 1.1}}  {\it Let $M$ be an
$n$-dimensional $(n\geq4)$ complete noncompact submanifold  with
flat normal bundle in sphere $S^{n+l}$. There exists a positive
constant $c(n)$ depending only on $n$ such that if the total
curvature is less than $c(n)$, then $H^p(L^2(M))=\{0\}$}, $2\leq p
\leq n-2$, where constant $c(n)$ is given by (3.8).

\vskip 5pt\noindent{\bf {Theorem 1.2}}  {\it Let $M$ be an
$n$-dimensional $(n\geq4)$ complete noncompact submanifold  with
flat normal bundle in sphere $S^{n+l}$. If the total curvature is
finite and $2\leq p \leq n-2$, then the dimension of $H^p(L^2(M))$
is finite.}

\begin{remark}
Theorem1.1 is a generalization of Theorem $B$. On the other hand,
harmonic $p$-forms $2\leq p \leq n-2$ are studied in Theorem 1.2
which is an extension of Theorem $A$. It is interesting to ask
whether there are finitely many non-parabolic ends on $M$ in Theorem
1.2.
\end{remark}

\section{Preliminaries}

Suppose $M$ is an $n$-dimensional complete  submanifold immersed in
an $n+l$ dimensional sphere $S^{n+l}$, $A$ is the second fundamental
form and $H$ is the mean curvature vector of $M$. The traceless
second fundamental form $\Phi$ is defined by
$$\Phi(X,Y)=A(X,Y)-\langle X,Y\rangle H,$$
for all vector field $X$ and $Y$, where $\langle,\rangle$ is the
metric of $M$. Obviously
$$|\Phi|^2=| A|^2-n|H|^2.$$
We say $M$ has finite total curvature if
$$\parallel\Phi\parallel_{L^n(M)}=(\int_{M}|\Phi|^n)^{\frac{1}{n}}<\infty.$$
$H^p(L^2(M))$ denotes the space of all $L^2$-harmonic $p$-forms on
$M$. Choose  local orthonormal frames $e_1,...,e_{n+l}$ on $S^{n+l}$
such that, restricted to $M$, $e_1,...,e_{n}$ are tangent to $M$.
Let $\omega_1,...,\omega_{n+l}$ are dual coframes. We then have
$\omega_\alpha=0$ for each $\alpha$, $n+1\leq\alpha\leq n+l$. From
Cartan's Lemma, we have $\omega_{\alpha i}=h^{\alpha}_{ij}\omega_j$.
$M$ has flat normal bundle implies that there exists an orthonormal
frame diagonalizing $h^\alpha_{ij}$ simultaneously.

We need the following Sobolev inequality which is a main tool in
proving our results.
 \vskip
5pt\noindent\begin{lemma}{(\cite{11,20})}\label{lem1} Let $M^n$ be
an n-dimensional complete noncompact oriented submanifold in sphere,
then
$$(\int_{M}| f|^{\frac{2n}{n-2}})^{\frac{n-2}{n}}\leq C_0[\int_{M}| \nabla f|^2+\int_{M}(| H|^2+1)f^2],$$
for all $f\in C^1_{0}(M)$, where $C_0$ depends only on $n$.
\end{lemma}

Lin, Han and Li proved the following estimate.
 \vskip
5pt\noindent\begin{lemma}{(\cite{13,22,23})}\label{lem1} Let $M^n$
be a complete submanifold with flat normal bundle in $S^{n+l}$,
$\omega$ be a $L^2$-harmonic $p$-form ($2\leq p\leq n-2$) on $M^n$,
then
$$|\omega|\Delta| \omega|\geq K_p|\nabla|\omega||^2+p(n-p)|\omega|^2+Q_p|\omega|^2,$$
where $Q_p={\inf}
_{i_1,...,i_n}(h^\alpha_{i_1i_1}+...+h^\alpha_{i_pi_p})(h^\alpha_{i_{p+1}i_{p+1}}+...+h^\alpha_{i_ni_n})$.
\end{lemma}

\section{Proof of our main Theorems }

{\bf {Proof of Theorem 1.1.}}
 $M$ has flat normal bundle implies
that there exists an orthonormal frame diagonalizing $h^\alpha_{ij}$
simultaneously. Choose proper local orthonormal frames,
$h^\alpha_{ij}$ are diagonalized simultaneously.
 Direct computation yields

\begin{align}
& 2\sum_{\alpha=n+1}^{n+l}(h^\alpha_{i_1i_1}+...+h^\alpha_{i_pi_p})(h^\alpha_{i_{p+1}i_{p+1}}+...+h^\alpha_{i_ni_n}) \notag  \\
  & = \sum_{\alpha=n+1}^{n+l}(h^\alpha_{i_1i_1}+...+h^\alpha_{i_ni_n})^2-\sum_{\alpha=n+1}^{n+l}(h^\alpha_{i_1i_1}+...+h^\alpha_{i_pi_p})^2
  -\sum_{\alpha=n+1}^{n+l}(h^\alpha_{i_{p+1}i_{p+1}}+...+h^\alpha_{i_ni_n})^2 \notag  \\
&\geq n^2| H|^2-\max\{p,n-p\}| A|^2  \notag  \\
&=\{p,n-p\}n| H|^2  -\max\{p,n-p\}| \Phi|^2 .
\end{align}
Substituting (3.1) into Lemma 2.2, we have
\begin{align}
|\omega|\Delta| \omega| \geq &
K_p|\nabla|\omega|\mid^2+p(n-p)|\omega|^2 \notag \\
&+\min\{p,n-p\}\frac{n}{2}| H\mid^2|\omega|^2
-\frac{1}{2}\max\{p,n-p\}| \Phi|^2|\omega|^2.
\end{align}
This together with the condition of $2\leq p\leq n-2$
 yields
\begin{equation}
|\omega|\Delta| \omega|\geq
\frac{1}{n-2}|\nabla|\omega||^2+2(n-2)|\omega|^2 +n| H|^2|\omega|^2
-\frac{n-2}{2}| \Phi|^2|\omega|^2.
\end{equation}
Setting $\eta \in C^\infty_0(M)$, multiplying (3.3) by $\eta^2$ and
integrating over $M$, we obtain

\begin{align}
\frac{n-2}{2}\int_{M}| \Phi|^2|\omega|^2\eta^2\geq  &
\frac{n-1}{n-2}\int_{M}|\nabla\mid\omega||^2\eta^2+2(n-2)\int_M|\omega|^2\eta^2 \notag \\
&+n\int_{M}| H|^2|\omega|^2\eta^2+2\int_{M}\eta|\omega|\langle\nabla
\eta,\nabla|\omega|\rangle.
\end{align}
Combining the H\"{o}lder inequality with Lemma 2.1, we get
\begin{align}
\int_{M}| \Phi|^2|\omega|^2\eta^2 \notag \leq & (\int_{M}|
\Phi|^n)^{\frac{2}{n}}(\int_{M}(|
\omega|\eta)^{\frac{2n}{n-2}})^{\frac{n-2}{n}} \notag \\
\leq & C_0(\int_{M}|
\Phi|^n)^{\frac{2}{n}}[\int_{M}|\nabla(\eta|\omega|)|^2+\int_{M}(|
H \mid^2+1)|\omega|^2\eta^2]\notag \\
\leq & C_0(\int_{M}|
\Phi|^n)^{\frac{2}{n}}[\int_{M}(|\nabla|\omega||^2\eta^2+|\omega|^2|\nabla\eta|^2+
2|\omega|\eta\langle\nabla\eta,\nabla|\omega|\rangle) \notag \\
&+\int_{M}(| H |^2+1)|\omega|^2\eta^2].
\end{align}
Setting $E=\frac{n-2}{2}C_0(\int_{M}| \Phi|^n)^{\frac{2}{n}}$ and
using (3.4) and (3.5) it follows that
\begin{align}
&E \int_{M}|
\omega|^2|\nabla\eta|^2+2(E-1)\int_M|\omega|\eta\langle\nabla\eta,\nabla|\omega|\rangle
\notag \\
& \geq
(\frac{n-1}{n-2}-E)\int_{M}|\nabla|\omega||^2\eta^2+[2(n-2)-E]\int_M|\omega|^2\eta^2
\notag \\
& +(n-E)\int_{M}| H|^2|\omega|^2\eta^2.
\end{align}
Using the Cauchy-Schwarz inequality in (3.6) , we get
\begin{align}
&(E+\frac{| E-1|}{\varepsilon }) \int_{M}|
\omega|^2|\nabla\eta|^2 \notag \\
& \geq (\frac{n-1}{n-2}-E-|
E-1|\varepsilon)\int_{M}|\nabla|\omega||^2\eta^2+[2(n-2)-E]\int_M|\omega|^2\eta^2
\notag \\
& +(n-E)\int_{M}| H|^2|\omega|^2\eta^2.
\end{align}
If
\begin{equation}
(\int_{M}|\Phi|^n)^{\frac{1}{n}}<\frac{2}{n-2}\sqrt{\frac{n-1}{2C_0}}=c(n),
\end{equation}
then
 $$\frac{n-1}{n-2}-E>0.$$  Choosing sufficient small $\varepsilon$,
we obtain
$$\frac{n-1}{n-2}-E-| E-1|\varepsilon>0, \ \  n-E>0,\ \
2(n-2)-E>0.$$  Let $\rho(x)$ be the geodesic distance on $M$ from
$x_0$ to $x$ and $B_r(x_0)=\{x\in M: \rho(x)\leq r\}$ for some fixed
point $x_0\in M$. Choose $\eta \in C^\infty _{0}(M)$ as
\begin{equation*}
 \eta=
  \begin{cases}
    1                              &\text{on}\ B_r(x_0),      \\
    0                              &\text{on}\ M\backslash B_{2r}(x_0), \\
    |\nabla\eta|\leq\frac2r\qquad  &\text{on}\ B_{2r}(x_0)\backslash B_{r}(x_0),
  \end{cases}
\end{equation*}
and $0\le\eta\le1$. Substituting the above $\eta$ into (3.7), we
finally have
\begin{align}
&\frac{4}{r^2}(E+\frac{| E-1|}{\varepsilon }) \int_{B_{2r}(x_0)}|
\omega|^2  \notag \\
& \geq (\frac{n-1}{n-2}-E-|
E-1|\varepsilon)\int_{B_r(x_0)}|\nabla|\omega||^2
+[2(n-2)-E]\int_{B_r(x_0)}|\omega|^2
\notag \\
& +(n-E)\int_{B_r(x_0)}| H|^2|\omega|^2. \notag
\end{align}
Since $\int_{M}| \omega|^2 <\infty$, by taking $r\rightarrow\infty$,
we have $\nabla|\omega |=0$ and $\omega=0$. That is
$H^p(L^2(M))=\{0\}$. This completes the proof of Theorem 1.1.

{\bf {Proof of Theorem 1.2.}} Let $\omega\in H^p(L^2(M)),\ \ 2\leq
p\leq n-2$ and $\eta \in C_0^{\infty}(M\backslash B_{r}(x_0)) $.
Similar to (3.7), we get
\begin{align}
&(F+\frac{\mid F-1\mid}{\varepsilon }) \int_{M\backslash
B_{r}(x_0)}|
\omega|^2|\nabla\eta|^2 \notag \\
& \geq (\frac{n-1}{n-2}-F-| F-1|\varepsilon)\int_{M\backslash
B_{r}(x_0)}|\nabla|\omega||^2\eta^2+[2(n-2)-F]\int_{M\backslash
B_{r}(x_0)}|\omega|^2\eta^2
\notag \\
& +(n-F)\int_{M\backslash B_{r}(x_0)}| H|^2|\omega|^2\eta^2,
\end{align}
where $F=\frac{n-2}{2}C_0(\int_{M\backslash B_{r}(x_0)}|
\Phi|^n)^{\frac{2}{n}}.$ The condition
$(\int_{M}|\Phi|^n)^{\frac{1}{n}}<\infty$ implies that there is a
decreasing positive function $\varepsilon(r)$ satisfying
$$\lim_{r\rightarrow\infty}\varepsilon(r)=0, \ \ (\int_{M\backslash B_{r}(x_0)}|\Phi|^n)^{\frac{2}{n}}<\varepsilon(r).$$
Thus we can choose $r=r_0>0$ such that
$$\frac{n-1}{n-2}-F=\frac{n-1}{n-2}-\frac{n-2}{2}C_0(\int_{M\backslash B_{r_0}(x_0)}\mid
\Phi\mid^n)^{\frac{2}{n}}>0.$$ Choosing  sufficient small
$\varepsilon$, we get
$$\frac{n-1}{n-2}-F-\mid
F-1\mid\varepsilon>0.$$ This together with (3.9) yields that
\begin{align}
\int_{M\backslash B_{r_0}(x_0)}|\nabla|\omega||^2\eta^2 &\leq
\frac{F+\frac{| F-1|}{\varepsilon }}{\frac{n-1}{n-2}-F-|
F-1|\varepsilon} \int_{M\backslash B_{r_0}(x_0)}|
\omega|^2|\nabla\eta|^2 \notag \\
&=C_1\int_{M\backslash B_{r_0}(x_0)}| \omega|^2|\nabla\eta|^2,
\end{align}

\begin{align}
\int_{M\backslash B_{r_0}(x_0)}|\omega|^2\eta^2 &\leq
\frac{F+\frac{| F-1|}{\varepsilon }}{2(n-2)-F} \int_{M\backslash
B_{r_0}(x_0)}|
\omega|^2|\nabla\eta|^2 \notag \\
&\leq C_1\int_{M\backslash B_{r_0}(x_0)}| \omega|^2|\nabla\eta|^2,
\end{align}

\begin{align}
\int_{M\backslash B_{r_0}(x_0)}| H|^2|\omega|^2\eta^2 &\leq
\frac{F+\frac{| F-1|}{\varepsilon }}{n-F} \int_{M\backslash
B_{r_0}(x_0)}|
\omega|^2|\nabla\eta|^2 \notag \\
&\leq C_1\int_{M\backslash B_{r_0}(x_0)}| \omega|^2|\nabla\eta|^2,
\end{align}
where the positive constant $C_1$ depends only on $n$. Applying
lamma 2.1 to $\eta|\omega|$
 and combining with (3.10),(3.11) and (3.12), we obtain
\begin{align}
&\int_{M\backslash B_{r_0}(x_0)}(\eta|\omega|)^{\frac{2n}{n-2}})^{\frac{n-2}{n}} \notag \\
&\leq C_0\int_{M\backslash
B_{r_0}(x_0)}[|\nabla|\omega||^2\eta^2+|\omega|^2|\nabla\eta|^2+
2|\omega|\eta\langle\nabla\eta,\nabla|\omega|\rangle+(| H
|^2+1)|\omega|^2\eta^2] \notag \\
&\leq C_0\int_{M\backslash
B_{r_0}(x_0)}[2|\nabla|\omega||^2\eta^2+2|\omega|^2|\nabla\eta|^2+(|
H |^2+1)|\omega|^2\eta^2] \notag \\
&\leq C_2\int_{M\backslash B_{r_0}(x_0)}|\omega|^2|\nabla\eta|^2,
\end{align}
where positive constant $C_2$ depends only on $n$.

Choose $\eta \in C^\infty _{0}(M\backslash B_{r_0}(x_0))$ as
\begin{equation*}
 \eta=
  \begin{cases}
    0                              &\text{on}\ B_{r_0}(x_0),      \\
    \rho(x)-r_0                    &\text{on}\ B_{r_0+1}(x_0)\backslash B_{r_0}(x_0), \\
    1                              &\text{on}\ B_{r}(x_0)\backslash B_{r_0+1}(x_0), \\
   \frac{2r-\rho(x)}{r}            &\text{on}\ B_{2r}(x_0)\backslash B_{r}(x_0), \\
    0                              &\text{on}\ M\backslash B_{2r}(x_0),
  \end{cases}
\end{equation*}
where  $\rho(x)$ is the geodesic distance on $M$ from $x_0$ to $x$
and $r>r_0+1$. By substituting $\eta$ into (3.13) it follows that
\begin{equation}
\int_{B_{r}(x_0)\backslash
B_{r_0+1}(x_0)}(|\omega|^{\frac{2n}{n-2}})^{\frac{n-2}{n}}\leq C_2
\int_{B_{r_0+1}(x_0)\backslash
B_{r_0}(x_0)}|\omega|^2+\frac{C_2}{r^2}\int_{B_{2r}(x_0)\backslash
B_{r}(x_0)}|\omega|^2.
\end{equation}
Since $|\omega|\in L^2(M),$ Letting $r\rightarrow \infty$, we
conclude that
\begin{equation}
\int_{B_{r}(x_0)\backslash
B_{r_0+1}(x_0)}(|\omega|^{\frac{2n}{n-2}})^{\frac{n-2}{n}}\leq C_2
\int_{B_{r_0+1}(x_0)\backslash B_{r_0}(x_0)}|\omega|^2.
\end{equation}
On the other hand , the H\"{o}lder inequality asserts that
\begin{equation}
\int_{B_{r_0+2}(x_0)\backslash B_{r_0+1}(x_0)}|\omega|^2\leq
(vol(B_{r_0+2}(x_0)))^\frac{2}{n}\int_{B_{r_{0}+2}(x_0)\backslash
B_{r_0+1}(x_0)}(|\omega|^{\frac{2n}{n-2}})^{\frac{n-2}{n}}
.\end{equation}
 From (3.15) and (3.16),  we conclude that there
exists a constant $C_3>0$ depending on $vol(B_{r_0+2}(x_0))$ and $n$
such that
\begin{equation}
\int_{B_{r_0+2}(x_0)}|\omega|^2\leq
C_3\int_{B_{r_0+1}(x_0)}|\omega|^2 .\end{equation} Fix a point $x
\in M$ and take $\tau \in C^1_0(B_1(x))$. Multiplying (3.3) by
$|\omega|^{q-2}\tau^2$ with $q>2$ and integrating by parts on
$B_{1}(x)$, we obtain
\begin{align}
&-2\int_{B_1(x)}\tau|\omega|^{q-1}\langle\nabla
\tau,\nabla|\omega|\rangle+\frac{n-2}{2}\int_{B_1(x)}| \Phi|^2|\omega|^q\tau^2 \notag \\
&\geq
(\frac{1}{n-2}+q-1)\int_{B_1(x)}|\omega|^{q-2}|\nabla|\omega||^2\tau^2
+2(n-2)\int_{B_1(x)}|\omega|^q\tau^2\notag \\
&+n\int_{B_1(x)}| H|^2|\omega|^q\tau^2
\end{align}
By using the Cauchy-Schwarz inequality it follows that
\begin{equation}
-2\int_{B_1(x)}\tau|\omega|^{q-1}\langle\nabla
\tau,\nabla|\omega|\rangle \leq
\frac{1}{n-2}\int_{B_1(x)}|\omega|^{q-2}|\nabla|\omega||^2\tau^2+(n-2)\int_{B_1(x)}|\omega|^{q}|\nabla\tau|^2.
\end{equation}
It follows from (3.18) and (3.19) that
\begin{align}
&(n-2)\int_{B_1(x)}|\omega|^{q}|\nabla\tau|^2 +\frac{n-2}{2}\int_{B_1(x)}| \Phi|^2|\omega|^q\tau^2 \notag \\
&\geq (q-1)\int_{B_1(x)}|\omega|^{q-2}|\nabla|\omega||^2\tau^2
+2(n-2)\int_{B_1(x)}|\omega|^q\tau^2\notag \\
&+n\int_{B_1(x)}| H|^2|\omega|^q\tau^2.
\end{align}
 On the other hand, setting  $f\in C^1_0(B_1(x))$, similar to
 Lemma 2.1, we have
\begin{equation}
(\int_{B_1(x)}| f|^{\frac{2n}{n-2}})^{\frac{n-2}{n}}\leq
C_0[\int_{B_1(x)}| \nabla f|^2+\int_{B_1(x)}(| H|^2+1)f^2].
\end{equation}
Applying (3.21) to $\tau |\omega|^{\frac{q}{2}}$ , we obtain
\begin{align}
(\int_{B_1(x)}(\tau^2
|\omega|^{q})^{\frac{n}{n-2}})^{\frac{n-2}{n}}&\leq
C_0\int_{B_1(x)}| \nabla (\tau
|\omega|^{\frac{q}{2}})|^2+C_0\int_{B_1(x)}(|
H|^2+1)\tau^2 |\omega|^{q}\notag \\
& \leq 2C_0\int_{B_1(x)}| \nabla \tau \mid^2|
\omega|^{q}+\frac{q^2}{2}C_0\int_{B_1(x)}\tau^2|\omega|^{q-2}|
\nabla
|\omega||^{2}\notag \\
& +C_0\int_{B_1(x)}(| H|^2+1)\tau^2 |\omega|^{q}.
\end{align}
Inequality (3.22) and (3.20) imply that
\begin{align}
&(\int_{B_1(x)}( \tau^2
|\omega|^{q})^{\frac{n}{n-2}})^{\frac{n-2}{n}}\notag \\
 &\leq
2C_0\int_{B_1(x)}|\nabla \tau
|^2|\omega|^{q}+\frac{q^2}{2(q-1)}C_0\int_{B_1(x)}[(n-2)|\nabla\tau|^2+\frac{n-2}{2}|\Phi|^2\tau^2]|\omega|^q \notag \\
&-\frac{q^2}{2(q-1)}C_0\int_{B_1(x)}[2(n-2)+n|
H|^2]|\omega|^{q}\tau^2 +C_0\int_{B_1(x)}(|
H|^2+1)\tau^2 |\omega|^{q} \notag \\
&\leq qC_4\int_{B_1(x)}(| \nabla \tau |^2+|\Phi|^2\tau^2)|\omega|^q,
\end{align}
where $C_4$ is a positive constant depending only on $n$. Let
$q_k=\frac{2n^k}{(n-2)^k}$ and $r_k=\frac{1}{2}+\frac{1}{2^{k+1}}$
for an integer $k\geq 0$. Choose $\tau_k\in C^\infty_0(B_{r_k}(x))$
such that $\tau_k=1$ on $B_{r_{k+1}}(x)$ and
$|\nabla\tau_k|\leq2^{k+3}$. Replacing $q$ and $\tau$ in (3.23) by
$q_k$ and $\tau_k$ respectively, we obtain
\begin{equation}
(\int_{B_{r_{k+1}}(x)}| \omega|^{q_{k+1}})^{\frac{1}{q_{k+1}}}\leq
[q_kC_4(4^{k+3}+\sup_{B_{1}(x)}|\Phi|^2)]^{\frac{1}{q_k}}(\int_{B_{r_k}(x)}|
\omega|^{q_k})^{\frac{1}{q_k}}.
\end{equation}
Apply the Morse interation to $\mid\omega\mid$ via (3.24), we
conclude that
$$\|\omega\|^2_{L^\infty(B_{\frac{1}{2}}(x))}\leq C_5
\int_{B_{1}(x)}|\omega|^2,$$ where $C_5$ is a positive constant
depending only on $n$ . Obviously
\begin{equation}|\omega(x)|^2\leq C_5
\int_{B_{1}(x)}|\omega|^2.\end{equation}
 Choose $x\in
\overline{B_{r_0+1}(x_0)}$ such that
$$|\omega(x)|^2=\|\omega\|^2_{L^\infty(B_{r_0+1}(x_0))}.$$
This together with (3.25) yields that
\begin{equation}
\|\omega\|^2_{L^\infty(B_{r_0+1}(x_0))}= |\omega(x)|^2\leq C_5
\int_{B_{1}(x)}|\omega|^2\leq C_5 \int_{B_{r_0+2}(x_0)}|\omega|^2.
\end{equation}
This together with (3.17) implies that there exists a positive
constant $C_6$ depending on $n$ and $vol(B_{r_0+2}(x_0))$, such that
\begin{equation}
\sup_{B_{r_0+1}(x_0)}|\omega|^2\leq C_6
\int_{B_{r_0+1}(x_0)}|\omega|^2.
\end{equation}
Let $\varphi$ be a finite dimensional subspace of $H^p(L^2(M))$.
 Lemma 11 in [25] implies that there exits $\omega\in \varphi$ such that
$$\frac{dim \varphi }{vol(B_{r_0+1}(x_0))}\int_{B_{r_0+1}(x_0)}|\omega|^2\leq |\{(^n_p), dim \varphi\}  \sup_{B_{r_0+1}(x_0)}|\omega|^2.$$
This together with (2.27) yields $dim \varphi\leq C_7$, where $C_7$
depending on $n$ and $vol(B_{r_0+1}(x_0))$. This implies that $dim
H^p(L^2(M))<\infty$, which completes the proof of Theorem 1.2.

%
\bibliographystyle{amsplain}

\end{document}